\documentstyle[11pt]{article}

\newtheorem{theorem}{Theorem}[section]
\newtheorem{lemma}[theorem]{Lemma}
\newtheorem{corollary}[theorem]{Corollary}
\newtheorem{proposition}[theorem]{Proposition}

\newtheorem{rem}[theorem]{Remark}
\newenvironment{remark}{\begin{rem}\rm}{\end{rem}}

\newtheorem{define}[theorem]{Definition}
\newenvironment{definition}{\begin{define}\rm}{\end{define}}

\newtheorem{ex}[theorem]{Example}

\newenvironment{proof}{\par {\em Proof.}}{ \hfill }

\font\twlmsbm=msbm10 scaled \magstep1
\font\egtmsbm=msbm8
\font\sixmsbm=msbm6
\newfam\msbmfam
 
\textfont\msbmfam=\twlmsbm
\scriptfont\msbmfam=\egtmsbm
\scriptscriptfont\msbmfam=\sixmsbm
 
\newcommand{\id}{\mbox{id}}

\newcommand{\Index}{\mbox{Index}\,}

\newcommand{\eps}{\varepsilon}
\newcommand{\lact}{\triangleright}

\newcommand{\la}{{<}}  
\newcommand{\ra}{{>}}  
\newcommand{\flip}{\varsigma}

\renewcommand{\star}{{\dag}}
\newcommand{\rtimes}{{>\!\!\!\triangleleft}}

\newcommand{\half}{\frac{1}{2}}

\newcommand{\1}{_{(1)}}
\newcommand{\2}{_{(2)}}
\newcommand{\3}{_{(3)}}

\newcommand{\tDelta}{\tilde\Delta}
\newcommand{\tS}{\tilde{S}}
\newcommand{\teps}{\tilde\eps}

\newcommand{\tI}{_{(\tilde{1})}}
\newcommand{\tII}{_{(\tilde{2})}}

\font\twlmsbm=msbm10 scaled \magstep1
\font\egtmsbm=msbm8
\font\sixmsbm=msbm6
\newfam\msbmfam
 
\textfont\msbmfam=\twlmsbm
\scriptfont\msbmfam=\egtmsbm
\scriptscriptfont\msbmfam=\sixmsbm
 
\def\Bbb#1{{\fam\msbmfam\relax#1}}  
 
\title{\bf A CHARACTERIZATION OF 
DEPTH 2 SUBFACTORS OF II${}_1~$FACTORS}
 
\author{Dmitri Nikshych\thanks{UCLA, Department of Mathematics, 
405 Hilgard Avenue, Los Angeles, CA
90095-1555; E-mail: nikshych@math.ucla.edu} \and Leonid Vainerman\thanks
{URA CNRS (Case 191), Universit\'e Pierre et Marie Curie, 4, place Jussieu,
F-75252, Paris Cedex 05 France; E-mail: vainerma@math.jussieu.fr}}
 
\date{October 15, 1998}

\begin{document}

\maketitle
\vskip 0.5cm
\begin{abstract}
We characterize finite index depth 2 inclusions of type II${}_1~$factors in
terms of actions of weak Kac algebras and weak $C^*$-Hopf algebras.
If $N\subset M \subset M_1 \subset M_2\subset \dots$ is the Jones tower 
constructed from such an inclusion $N\subset M$, then $B=M^\prime \cap M_2$ 
has a natural structure of a weak $C^*$-Hopf algebra and 
there is a minimal action of $B$ on $M_1$ such that $M$ is the fixed point 
subalgebra of $M_1$ and $M_2$ is isomorphic to the crossed product of $M_1$ 
and $B$. This extends the well-known results for irreducible depth 2 inclusions.
\end{abstract}


\begin{section}
{Introduction}

Let $N\subset M$ be a finite index depth 2 inclusion of type II${}_1$
factors  and $N\subset M \subset M_1 \subset M_2 \subset \dots$ the   
corresponding Jones tower.
It was announced by A.~Ocneanu and was proved in \cite{S}, \cite{D}, \cite{L}
that if $N\subset M$ is irreducible, i.e., such that $N'\cap M=\Bbb{C}$, then
$B=M^\prime\cap M_2$ has a natural structure of a finite-dimensional
Kac algebra and there is a canonical  outer action of $B$ on $M_1$  
such that $M=M_1^B$, the fixed point subalgebra of $M_1$ with respect to this
action, and $M_2$ is isomorphic to the crossed product $M_1 \rtimes B$.
The outerness condition is equivalent to the relative commutant 
$M_1^\prime \cap M_1 \rtimes B$ being trivial (such actions are also called 
minimal). In the case of an infinite index a similar description in terms of
multiplicative unitaries and quantum groups was obtained in \cite{EN}.
 
In this work we extend the above result to (in general, reducible, i.e.,
such that $\Bbb{C}\subset N'\cap M$) finite index depth 2 inclusions of type
II${}_1~$factors. We replace usual Kac algebras (Hopf $C^*$-algebras)
by weak Kac algebras \cite{NV} or weak $C^*$-Hopf algebras \cite{BNSz}. 
A weak Kac algebra is a special case of a weak $C^*$-Hopf algebra
characterized by the property $S^2=\id$. Weak Kac algebras naturally
arise in the situations when the index $[M:N]$ is integer, e.g., when
the inclusion is given by the crossed product with a usual Kac algebra.
It was shown in \cite{NV} that the
category of weak Kac algebras is equivalent to those
of generalized Kac algebras of T.~Yamanouchi \cite{Y} (another proof
of that can be found in \cite{Nill}) and of Kac bimodules
(an algebraic version of Hopf bimodules of J.-M. Vallin
\cite{Val}). The advantage 
of the language of weak Kac algebras and weak $C^*$-Hopf
algebras is that their defining  axioms are clearly self-dual, so
it is easy to work with both weak Kac algebra (weak $C^*$-Hopf algebra) and 
its dual simultaneously.

Let us mention that a possibility of characterizing finite index depth 
2 inclusions in terms of weak $C^*$-Hopf algebras was suggested in \cite{NSzW}. 
For an arbitrary (possibly infinite) index M.~Enock and J.-M.~Vallin have
obtained a similar description in terms of pseudo-multiplicative unitaries
\cite{EVal}.
\medskip
 
The paper is organized as follows.
 
In Section 2 (Preliminaries) we briefly review, following \cite{NV},
\cite{BNSz} and \cite{NSzW}, the basic definitions and facts of the theory 
of weak Kac algebras and weak $C^*$-Hopf algebras, including their actions on
von Neumann algebras.
 
Section 3 is devoted to
establishing a non-degenerate duality between the finite dimensional
$C^*$-algebras $A=N'\cap M_1$ and $B=M'\cap M_2$, which
gives a natural coalgebra structures on them.

In Sections 4 and 5 we investigate the relations between algebra and coalgera
structures on $B$, following the general strategy of Szymanski's reasoning
\cite{S} based on the above duality. It turns out that the square of the
corresponding antipode is implemented by a positive invertible element
determined by $\Index\,\tau\vert_{M^\prime\cap M_1}$, the Watatani index 
\cite{W} of the restriction of the Markov trace $\tau$ on $M^\prime\cap M_1$.
That is why it is natural to consider the cases of scalar and non-scalar
$\Index\,\tau\vert_{M^\prime\cap M_1}$ in which the antipode is respectively
involutive and non-involutive. The main result here is that in the mentioned
cases $B$ and $A$ are biconnected weak Kac algebras and weak $C^*$-Hopf 
algebras respectively (they are usual Kac algebras iff the inclusion 
$N\subset M$ is irreducible). 
We also prove in Section 4, that if $[M:N]$ is an integer which
has no divisors of the form $n^2,\ n>1$, then the inclusion is irreducible and
$B$ is a Kac algebra acting outerly on $M_1$. In particular, if $[M:N]=p$ is
prime, then $B$ must be the group algebra of the cyclic group $G=Z/pZ$.
 
Finally, in Section  6 we show that there exists a canonical
(left) minimal action of $B$ on $M_1$ such that $M$ is the fixed point
subalgebra of $M_1$ with respect to this action, and  $M_2$ is isomorphic
to $M_1\rtimes B$, the crossed product of $M_1$ and $B$. 
The minimality condition means that  the relative commutant 
$M_1^\prime \cap M_1 \rtimes B$ is minimal possible, in which case
it is isomorphic to the Cartan subalgebra $B_s\subset B$. 

It is important to stress that in the above situation one can take
$$
\begin{array}{ccc}
B^* & \subset & B^* \rtimes B  \\
\cup&         & \cup        \\
B^*\cap B & \subset & B,
\end{array}
$$
where $B^*=A$, as a canonical commuting square \cite{Popa2} of the inclusion
$M_1 \subset M_2$. The above square, and thus the equivalence class 
of inclusions, is completely determined by $B$.
This implies that every biconnected weak $C^*$-Hopf algebra has at most
one minimal action on a given II$_{1}$ factor and thus correspond to no
more than one (up to equivalence) finite index depth 2 subfactor.
Note that any  biconnected weak Kac algebra admits a unique minimal action
on the hyperfinite II${}_1~$factor \cite{Nik}.

Let us remark that this characterization of depth 2
inclusions means that weak Kac algebras provide a good setting
for studying actions of usual Kac algebras on II${}_1$~factors,
since any (not necessarily minimal) action of a Kac algebra produces
a depth 2 inclusion and one can canonically associate 
with this action a weak Kac algebra completely describing it. 
More details on this will be published elsewhere.
 
{\bf Acknowledgements.} The first author would like to thank S.~Popa
and E.~Effros for helpful advices, and E.~Vaysleb for his useful
comments on the preliminary version of this paper.
The second author is deeply grateful to M.~Enock and J.-M.~Vallin for
many important discussions.
\end{section}
 

\begin{section}
{Preliminaries}

Our main references to
finite dimensional weak $C^*$-Hopf algebras are \cite{BNSz} and \cite{Nill}. 
Weak Kac algebras, a special case of this 
notion characterized by the property $S^2=\id$, were considered in \cite{NV}.
These objects generalize both finite groupoid algebras and usual Kac algebras.

A {\em weak Kac algebra} $B$ is a finite dimensional $C^*$-algebra equipped
with the {\em comultiplication} $\Delta : B\to B \otimes B$,
{\em counit} $\eps : B \to \Bbb{C}$, and {\em antipode} $S : B\to B$,
such that ($\Delta,\,\eps$) defines a coalgebra structure on $B$ and
the following axioms hold for all $b,\,c\in B$ (we use Sweedler's notation
$\Delta(b) =b\1\otimes b\2$ for the comultiplication) :

\begin{enumerate}
\item[(1)] $\Delta$ is a $*$-preserving (but not necessarily unital)
homomorphism : 
$$
\Delta(bc) = \Delta(b)\Delta(c), \quad \Delta(b^*) = \Delta(b)^{*\otimes *},
$$
\item[(2)] The target counital map $\eps^t$, defined by 
$\eps^t(b)=\eps(1\1b)1\2$, satisfies the relations
$$
b\eps^t(c) = \eps(b\1c)b\2, \qquad
b\1 \otimes \eps^t(b\2)  = 1\1b\otimes 1\2,
$$
\item[(3)] 
$S$ is an anti-algebra  and anti-coalgebra  map
such that $S^2=\id$, $(S\circ *) =(*\circ S)$, and
$$
b\1 S(b\2) = \eps^t(b).
$$
\end{enumerate}
\medskip

If instead of the conditions  $S^2=\id$ and  $(S\circ *) =(*\circ S)$
we have a less restrictive property $(S\circ *)^2 = \id$, then $B$
is called a {\em weak $C^*$-Hopf algebra}.
\medskip

Note that the axioms (2) and (3) above are equivalent to the following axioms
for the source counital map $\eps^s(b)=1\1\eps(b1\2)$ :
\begin{enumerate}
\item[(2')] 
$ \eps^s(c)b = b\1\eps(cb\2), \qquad
\eps^s(b\1) \otimes b\2 = 1\1 \otimes b1\2, $
\item[(3')] 
$S(b\1)b\2 = \eps^s(b)$.
\end{enumerate}

The dual vector space $B^*$ has a natural structure of a weak Kac algebra
(weak $C^*$-Hopf algebra) given by dualizing the structure operations of 
$B$, see \cite{BNSz}, \cite{NV}.

The main difference between weak Kac ($C^*$-Hopf) algebras and classical
Kac algebras is that the images of the counital maps are, in general,
non-trivial unital $C^*$-subalgebras of $B$, called {\em Cartan subalgebras} 
(note that we have $\eps^t\circ\eps^t = \eps^t$ and 
$\eps^s\circ\eps^s = \eps^s$) :
\begin{eqnarray*}
B_t &=& \{x\in B \mid \eps^t(x) =x \}
  =  \{x\in B \mid \Delta(x) = x 1\1 \otimes 1\2 = 1\1 x\otimes 1\2 \}, \\
B_s &=& \{x\in B \mid \eps^s(x) =x \} 
= \{x\in B \mid \Delta(x) = 1\1 \otimes x 1\2 = 1\1 \otimes 1\2 x \}.
\end{eqnarray*}
The Cartan subalgebras commute : $[B_t,\, B_s] =0$, also we have 
$S\circ\eps^s = \eps^t\circ S$ and $S(B_t) =B_s$. We say that $B$ is 
{\em connected} \cite{Nik} if $B_t \cap Z(B)= \Bbb{C}$ (where $Z(B)$
denotes the center of $B$), i.e.,\ if the inclusion $B_t \subset B$
is connected. $B$ is connected iff $B_t^* \cap B_s^* = \Bbb{C}$
(\cite{Nik}, Proposition 3.11). We say that $B$ is {\em biconnected} 
if both $B$ and $B^*$ are connected.
\medskip

Weak Kac ($C^*$-Hopf) algebras have integrals in the following sense.

\begin{enumerate}
\item[]
There exists a unique projection $p\in B$,
called a {\em Haar projection}, such that for all $x\in B$ :
$$
x p  = \eps^t(x)p ,\quad S(p) = p, \quad \eps^t(p)=1.
$$
\item[]
There exists a unique positive functional $\phi$ on $B$,
called a {\em normalized Haar functional} (which is a trace iff
$B$ is a weak Kac algebra), such that
$$
(\id\otimes \phi)\Delta = (\eps^t\otimes \phi)\Delta,\quad
\phi\circ S =S,\quad
\phi\circ \eps^t = \eps.
$$
\end{enumerate}

The following notions of action, crossed product, and fixed point subalgebra
were introduced in \cite{NSzW}.

A (left) {\em action} of a weak Kac ($C^*$-Hopf algebra) $B$ on a
von Neumann algebra $M$ is a linear map
$$
B\otimes M \ni b\otimes x \mapsto (b\lact x)\in M
$$
defining a structure of a left $B$-module on $M$ such that for all
$b\in B$ the map $b\otimes x \mapsto (b\lact x)$ is weakly continuous and
\begin{enumerate}
\item[(1)]
$b\lact xy = (b\1\lact x)(b\2\lact y),$
\item[(2)]
$(b\lact x)^* = S(b)^* \lact x^*, $
\item[(3)]
$b\lact 1 = \eps^t(b)\lact 1$, and  $b\lact 1 =0$ iff $\eps^t(b)=0.$ 
\end{enumerate}
\medskip

A {\em crossed product} algebra 
$M\rtimes B$ is constructed as follows. As a $\Bbb{C}$-vector space it is 
$M \otimes_{B_t} B$, where  $B$ is a left $B_t$-module  
via multiplication and $M$ is a right $B_t$-module via multiplication 
by the image of $B_t$ under $z\mapsto (z\lact 1)$; that is, we identify 
$$
x(z\lact 1) \otimes b  \equiv  x \otimes zb  
$$
for all $x\in M,\,b\in B,\,z\in B_t$.
Let $[x\otimes b]$  denote the class of $x\otimes b$. 
A $*$-algebra structure on $M\rtimes B$ is defined by
$$
[x\otimes b][y\otimes c] = [x(b\1\lact y) \otimes b\2 c], \quad
[x\otimes b]^* = [(b\1^*\lact x^*) \otimes b\2^*] \\
$$
for all $x,y\in A,\,b,c \in B$. It is possible to show that this
abstractly defined $*$-algebra $M\rtimes B$ is *-isomorphic to a 
weakly closed algebra of operators on some Hilbert space \cite{NSzW},
i.e., $M\rtimes B$ is a von Neumann algebra.

The collection
$ M^B =\{ x\in M \mid b\lact x = \eps^t(b)\lact x,\ \forall b\in B\}$
is a von Neumann subalgebra of $M$, called a
{\em fixed point subalgebra}.

The relative commutant $M^\prime \cap M \rtimes B$ always contains
a *-subalgebra isomorphic to $B_s$. Indeed, if $z\in B_s$, then
it follows easily from the axioms of a weak $C^*$-Hopf algebra
that $\Delta(z) =1\1 \otimes 1\2 z$, therefore
\begin{eqnarray*}
[1\otimes z][x \otimes 1] 
&=& [(z\1\lact x)\otimes z\2]  = [(1\1\lact x)\otimes 1\2z] \\
&=&  [x\otimes z] = [x \otimes 1][1\otimes z],
\end{eqnarray*}
for all $x\in M$, and $B_s \subset M^\prime \cap M\rtimes B$.
We say that the action $\lact$ is {\em minimal} if 
$B_s = M^\prime \cap M\rtimes B$. 

\end{section}


\begin{section}
{Duality between relative commutants}

Let $N \subset M$ be a depth 2 inclusion of type II${}_1~$factors
with a finite index $[M:N] =\lambda^{-1}$ and 
$$
N \subset M \subset M_1 \subset M_2 \subset \cdots
$$
be the corresponding Jones tower, 
$M_1 = \langle M,\,e_1 \rangle,\, M_2 = \langle M_1,\,e_2 \rangle,\dots,$ 
where $e_1 \in N^\prime \cap M_1,\,e_2 \in M^\prime \cap M_2,\cdots$ 
are the Jones  projections. The depth 2 condition means that 
$N^\prime \cap M_{2}$ is the basic construction of the inclusion 
$N^\prime \cap M  \subset  N^\prime \cap M_{1}$. Let $\tau$ be
the normalized (Markov) trace on $M_2$.

With respect to this trace,
the square of algebras in the upper right corner of the diagram below
$$
\begin{array}{ccccc}
N^\prime \cap M & \subset & N^\prime \cap M_{1} &
\subset &  N^\prime \cap M_{2} \\
   &   & \cup   &   & \cup \\
                    &         & M^\prime \cap M_{1} &
\subset &  M^\prime \cap M_{2} \\
   &   &    &   & \cup \\
   &  &  &  &  M_1^\prime \cap M_{2}.
\end{array}
$$
is commuting ($E_{M_1}\circ E_{M^\prime} =  E_{M^\prime} \circ E_{M_1}$
on $N^\prime \cap M_2$) and non-degenerate, i.e., 
$N^\prime \cap M_{2} =(N^\prime \cap M_1)(M^\prime \cap M_{2})$. 
This square is called a standard (or canonical) 
commuting  square of the inclusion $M_1\subset M_2$ \cite{Popa2}. 
 
Let us denote
\begin{eqnarray*}
A = N^\prime \cap M_1, &\qquad&  B = M^\prime \cap M_2, \\
A_t = N^\prime \cap M, \qquad  A_s = &M^\prime \cap M_1& =B_t,
\qquad B_s =  M_1^\prime \cap M_{2}.
\end{eqnarray*}
Note that 
$A_t$ commutes with $B$, $B_s$ commutes with $A$, and $A\cap B =A_s =B_t$.

The next lemma will be frequently used in the sequel without
specific reference.

\begin{lemma}
\label{pimsner-popa}
$(N^\prime \cap M_2)e_2 =Ae_2$ and $(N^\prime \cap M_2)e_{1} =Be_1$.
More precisely, for any $x \in N^\prime \cap M_2$ we have
$$
xe_2 = \lambda^{-1}E_{M_1}(xe_2)e_2, \qquad
xe_1 = \lambda^{-1}E_{M^\prime}(xe_1)e_1.
$$
\end{lemma}
\begin{proof}
This statement is a special case of (\cite{PP}, Lemma~1.2) since 
$N^\prime \cap M_2$ is the basic construction for the inclusions
$N^\prime \cap M \subset N^\prime \cap M_1$ 
and $M_1^\prime \cap M_2 \subset M^\prime \cap M_2$ with
the corresponding Jones projections $e_2$ and $e_1$ respectively.
\end{proof}  
\medskip

Let us denote $d =\dim( M^\prime \cap M_1)$.

\begin{proposition}
\label{duality}
The form 
$$
\la a,\, b\ra =  d\lambda^{-2}\tau(a e_2e_1b),\qquad
a\in A,\, b\in B
$$
defines a non-degenerate duality between $A$ and $B$.
\end{proposition}
\begin{proof}
If $a\in A$ is such that $\la a,\, B\ra = 0$ , then
$$
\tau(a e_2 e_1 B) = \tau(a e_2 e_1 (N^\prime \cap M_2)) = 0,
$$
therefore, using the Markov property of $\tau$ and properties of
Jones projections, we get
\begin{eqnarray*}
\tau(a a^*) &=& \lambda^{-1} \tau(a e_2 a^*) =
\lambda^{-2}\tau(a e_2 e_1 (e_2 a^* ) ) = 0,
\end{eqnarray*}
so $a=0$. Similarly for $b\in B$.
\end{proof}

\begin{definition}
Using the form $\la , \ra$ define the  
comultiplication $\Delta_B$, counit $\eps_B$,
and antipode $ S_B$ as follows :
\begin{eqnarray*}
\Delta_B : B \to B\otimes B :\quad
& & \la a_1a_2,\, b \ra =  \la a_1,\, b\1\ra \la a_2,\, b\2\ra, \\
\eps_B : B \to \Bbb{C} : \quad
& & \eps_B(b) = \la 1,\, b\ra = \lambda^{-1} d \tau(be_2), \\ 
S_B: B \to B  :\quad
& & \la a,\, S_B(b)\ra = \overline{\la a^*,\, b^* \ra},
\end{eqnarray*}
for all $a,a_1,a_2\in A$ and $b\in B$.
Similarly, we define $\Delta_A$, $\eps_A$, and  $S_A$.
\end{definition}
\medskip

Clearly, $(B,\,\Delta_B,\,\eps_B)$  (resp.\ $(A,\,\Delta_A,\,\eps_A)$)
becomes  a coalgebra. Let us investigate the relations between the 
algebra and coalgebra structures on $B$.

\end{section}


\begin{section}
{Weak Kac algebra structure on $M'\cap M_2$ 
(the case of a scalar Watatani index of $\tau\vert_{M^\prime\cap M_1}$)}

\begin{lemma}
\label{trick}
For all $a\in A$ and $b_1,b_2\in B$ we have
$$
\la a,\,b_1b_2 \ra =  \lambda^{-1} \la E_{M_1}(b_2ae_2),\, b_1\ra.
$$
\end{lemma}
\begin{proof}
Using the definition of $\la , \ra$ we have
\begin{eqnarray*}
\la a,\,b_1b_2 \ra &=& d\lambda^{-2}\tau(b_2a e_2 e_1 b_1) =
d\lambda^{-3} \tau(E_{M_1}(b_2ae_2) e_2 e_1  b_1 ) \\
&=& \lambda^{-1} \la E_{M_1}(b_2ae_2),\,b_1 \ra.
\end{eqnarray*}
\end{proof}

\begin{proposition}
\label{counital maps}
Let $\eps_B^t(b)  = \eps_B(1\1b)1\2$. Then
$\eps_B^t(b) = \lambda^{-1}E_{M_1}(be_2)$ and
$$
\la a,\, \eps_B^t(b) \ra = d\lambda^{-2}\tau(ae_1be_2) =
\lambda^{-1} \la E_M(ae_1),\,b \ra. 
$$
\end{proposition}
\begin{proof}
Using Lemma~\ref{trick}, definitions of $\Delta_B$ and $\eps_B$, we have
\begin{eqnarray*}
\la a,\, \eps_B(1\1b)1\2 \ra
&=& \la 1,\,1\1b \ra \la a,\, 1\2 \ra \\
&=& \la \lambda^{-1}E_{M_1}(be_2), 1\1 \ra \la  a,\, 1\2 \ra \\
&=& \la \lambda^{-1}E_{M_1}(be_2)a, 1 \ra 
 = \la a,\, \lambda^{-1}E_{M_1}(be_2) \ra, \\
\end{eqnarray*}
from where the first statement follows. For the second one, we have,
using the $\lambda$-Markov property and the fact that $e_2$ commutes 
with $M$,
\begin{eqnarray*}
\la a,\, \lambda^{-1}E_{M_1}(be_2) \ra 
&=& d \lambda^{-2} \tau( a e_1  e_2 \lambda^{-1}E_{M_1}(be_2)) \\
&=& d \lambda^{-2} \tau( a e_1  \lambda^{-1}E_{M_1}(be_2) e_2) \\
&=& d \lambda^{-2} \tau(a e_1 b e_2) 
 = d \lambda^{-3} \tau(E_M(a e_1)e_1 b e_2) \\
&=& d \lambda^{-3} \tau(E_M(a e_1) e_2e_1 b) 
 = \lambda^{-1} \la E_M(ae_1),\,b \ra. 
\end{eqnarray*}
\end{proof}

\begin{proposition}
\label{axiom A2}
For all $b,c\in B$ we have
$$ 
b\1 \otimes \eps_B^t(b\2) = 1\1b \otimes 1\2, \qquad
b\eps_B^t(c) = \eps_B(b\1 c)b\2.
$$
\end{proposition}
\begin{proof}
For all $a_1,\, a_2 \in A$ we compute, using Lemma~\ref{trick}
and Proposition~\ref{counital maps} :    
\begin{eqnarray*}
\la a_1,\, b\1 \ra \la a_2,\, \eps_B^t(b\2) \ra
&=& \la a_1 \lambda^{-1}E_M(a_2e_1),\, b\ra \\
&=& \lambda^{-2}  \la E_{M_1}(ba_1 E_M(a_2e_1)e_2),\, 1\ra \\
&=& \lambda^{-2}  \la E_{M_1}(ba_1e_2)E_M(a_2e_1),\, 1\ra \\
&=& d\lambda^{-3} \tau( E_{M_1}(ba_1e_2)E_M(a_2e_1)e_1) \\
&=& d\lambda^{-2} \tau( E_{M_1}(ba_1e_2)a_2e_1) \\
&=& \lambda^{-1} \la E_{M_1}(ba_1e_2)a_2, 1 \ra \\
&=& \la \lambda^{-1}E_{M_1}(ba_1e_2),\, 1\1 \ra \la a_2,\, 1\2 \ra \\
&=& \la a_1,\, 1\1b \ra \la a_2,\, 1\2 \ra, \\
\la  a,\, b\eps_B^t(c) \ra 
&=& \la \eps_B^t(c) a,\, b\ra \\
&=& \la \lambda^{-1}E_{M_1}(ce_2)a,\,b \ra \\
&=& \la \lambda^{-1}E_{M_1}(ce_2),\, b\1 \ra \la a,\, b\2 \ra \\
&=& \la 1,\, b\1c \ra \la a,\, b\2 \ra \\ 
&=& \la a,\, \eps_B(b\1 c)b\2 \ra.
\end{eqnarray*}
Since the duality is non-degenerate, the result follows.
\end{proof}
\medskip

The antipode map assigns to each $b\in B$
a unique element $S_B(b)\in B$ such that
$\tau(ae_2 e_1 S_B(b)) = \tau (b e_1 e_2 a)$ for all $a\in A$,
or, equivalently,
$$
E_{M_1}( b e_1 e_2) = E_{M_1}(e_2 e_1 S_B(b)).
$$
Taking $a=e_1$ and using the $\lambda$-Markov property of $e_1$
we get $\tau\circ S_B =\tau$. Similarly,
$E_{M^\prime}(S_A(a)e_2 e_1) = E_{M^\prime}(e_1 e_2 a)$
and $\tau\circ S_A =\tau$.

\begin{remark}
\label{even more}
Note that the condition 
$E_{M_1}( b e_1 e_2) = E_{M_1}(e_2 e_1 S_B(b))$
implies that
$$
E_{M_1}( b x e_2) = E_{M_1}(e_2 x S_B(b))\quad \mbox{for all } x\in M_1.
$$
Indeed, any $x\in M_1$ can be written as 
$x = \sum\, x_ie_1 y_i$ with $x_i,y_i \in M \subset B^\prime$.
Similarly, we have 
$$
E_{M^\prime}(S_A(a)ye_1) = E_{M^\prime}(e_1 ya) \quad
\mbox{for all } y\in M^\prime.
$$
\end{remark}

\begin{proposition}
\label{properties of S}
The following identities hold :
\begin{enumerate}
\item[(i)]
$S_B(b) = \lambda^{-3} E_{M^\prime}(e_1e_2 E_{M_1}(be_1e_2) ), $
\item[(ii)]
$S_B(B_s)=B_t$,
\item[(iii)]
$S_B^2(b) =b$ and $S_B(b)^* = S_B(b^*)$,
\item[(iv)]
$S_B(bc) = S_B(c)S_B(b)$ and 
$\Delta_B(S_B(b)) = \flip (S_B\otimes S_B)\Delta_B(b)$.
\end{enumerate}
\end{proposition}
\begin{proof}
(i) We have
\begin{eqnarray*}
S_B(b) 
&=& \lambda^{-1}  E_{M^\prime}( e_1 S_B(b)) 
    = \lambda^{-2} E_{M^\prime}( e_1 e_2 e_1 S_B(b)) \\
&=& \lambda^{-3} E_{M^\prime} ( e_1 e_2 E_{M_1}(e_2 e_1 S_B(b) )
    = \lambda^{-3} E_{M^\prime} ( e_1 e_2 E_{M_1}(b e_1 e_2) ).
\end{eqnarray*}
(ii) If $z\in B_s$ then $ze_2=e_2z$ and by
the explicit formula (i) we get,
\begin{eqnarray*}
S_B(z)
&=& \lambda^{-3} E_{M^\prime} ( e_1 e_2 E_{M_1}(e_1 z e_2) )
 =  \lambda^{-2} E_{M^\prime} ( e_1 E_{M_1}(z e_2) ) \\
&=&  \lambda^{-1} E_{M_1}(z e_2) = \eps_B^t(z) \in B_t.
\end{eqnarray*}
(iii) Since $E_{M_1}$ preserves $*$, we get
$E_{M_1}( e_2 e_1 b^*) = E_{M_1}(S_B(b)^* e_2 e_1),$ 
from where $S_B(S_B(b)^*)^* = b$. Next, using Lemma~\ref{trick},
Remark~\ref{even more}, and the $\lambda$-Markov property of
$e_2$, we compute
\begin{eqnarray*}
\tau(a e_2 e_1 b) 
&=&  \lambda^{-1} \tau(E_{M_1}(b a e_2)  e_2 e_1 ) 
      =  \lambda^{-1} \tau(E_{M_1}(e_2 a S_B(b))  e_2 e_1 ) \\ 
&=&  \lambda^{-1} \tau(e_2 E_{M_1}(e_2 a S_B(b))  e_1 ) 
      =  \tau( e_2 a S_B(b) e_1 ) \\ 
&=&  \tau( S_B(b) e_1 e_2 a) =   \tau(a e_2 e_1 S_B^2(b)).
\end{eqnarray*}
therefore, $S_B^2(b) =b$ and $S_B(b)^* = S_B(b^*)$.
\newline
(iv) Using Remark~\ref{even more}, we have
\begin{eqnarray*}
\tau( a e_2 e_1 S_B(bc))
&=& \tau( bc e_1 e_2 a) 
= \lambda^{-1} \tau(  c e_1 e_2 E_{M_1}(e_2ab)) \\
&=& \lambda^{-1} \tau(E_{M_1}(e_2ab)  e_2 e_1 S_B(c) )\\
&=& \lambda^{-1} \tau( E_{M_1}( S_B(b)a e_2)  e_2 e_1 S_B(c) )\\
&=& \tau( a e_2 e_1 S_B(c) S_B(b)),\\
\end{eqnarray*}
which proves that $\la a,\, S_B(bc) \ra = \la a,\, S_B(c) S_B(b) \ra$.
Similarly, one can prove that $S_A$ is anti-multiplicative,
and since $\la a,\, S_B(b) \ra = \la S_A(a),\, b\ra$,
the second part of (iv) follows.
\end{proof}
\medskip


Let $\{ f_{kl}^{\alpha} \}$ be a system of matrix
units in  $B_t = M^\prime \cap M_1 = \oplus_\alpha\, M_{m_\alpha}(\Bbb{C})$,
where $\sum m_\alpha^2 =d$, and let $\tau_\alpha = \tau(f_{kk}^{\alpha})$.

\begin{proposition}
\label{Delta 1} 
The explicit formula for $\Delta_B(1)$ is
$$\Delta_B(1) = \sum_{\alpha kl}\, \frac{1}{d\tau_\alpha}\,
S_B(f_{kl}^{\alpha}) \otimes f_{lk}^{\alpha}.
$$
In particular, $\Delta_B(1)$ is a positive element in $B_s\otimes B_t$.
\end{proposition}
\begin{proof}
Note that the map $x \mapsto  \sum_{\alpha kl}\,
\frac{\tau( x f_{lk}^{\alpha} ) }{\tau_\alpha} f_{kl}^{\alpha}$ 
defines the $\tau$-preserving conditional expectation on $B_t$.
For all $a_1, a_2 \in A$ we have
\begin{eqnarray*}
\lefteqn{ \sum_{\alpha kl}\, \frac{1}{d \tau_\alpha}
\la a_1,\, S_B(f_{kl}^{\alpha})\ra \la a_2,\, f_{lk}^{\alpha} \ra = } \\
&=& d^2\lambda^{-4}   \sum_{\alpha kl}\, \frac{1}{d \tau_\alpha}
\tau(a_1 e_2 e_1 S_B(f_{kl}^{\alpha})) \tau( a_2 e_2 e_1 f_{lk}^{\alpha}) \\
&=& d \lambda^{-3} \sum_{\alpha kl}\, 
\tau(f_{kl}^{\alpha} e_1e_2 a_1) 
\frac{\tau(a_2 e_1 f_{lk}^{\alpha}) }{\tau_\alpha}  \\
&=& d \lambda^{-3} \tau( E_{M^\prime}(a_2e_1) e_1 e_2 a_1 ) \\
&=& d \lambda^{-2}\tau( a_1 a_2 e_1 e_2 )  = \la a_1a_2,\, 1 \ra,
\end{eqnarray*}  
which proves the statement.
\end{proof}

\begin{corollary}
\label{s otimes id e}
$  \Delta_B(1) = \sum_{\alpha kl}\, \frac{1}{m_\alpha}\,
S_B(f_{kl}^{\alpha}) \otimes f_{lk}^{\alpha}H,$
where $H$ is canonically defined by
$$
H = S_B(1\1)1\2 = \frac{1}{d}\,\
\sum\nolimits_\alpha\, \frac{m_\alpha}{\tau_\alpha}\,
\sum\nolimits_k\,f_{kk}^{\alpha} 
= \frac{1}{d}\,\Index\tau\vert_{M^\prime\cap M_1} \in Z(B_t),
$$
where $\Index\tau\vert_{M^\prime\cap M_1}$ is the Watatani index \cite{W}
of the restriction of $\tau$ to $M^\prime\cap M_1$ and 
$Z(\cdot)$ denotes the center of the algebra. We also have $\tau(H)=1$.
\end{corollary}

\begin{proposition}
\label{my pain}
For all $b\in B$ we have 
$\eps_B^t(b\1)b\2 =Hb$.
\end{proposition}
\begin{proof}
Applying $E_{M^\prime}$ to both sides of
$E_{M_1}(b^* e_1 e_2) = E_{M_1}(e_2 e_1 S_B(b^*))$
and using the relation  
$E_{M_1}\circ E_{M^\prime} =  E_{M^\prime}\circ E_{M_1}$, we get
$$
E_{M_1}( b^* e_2) = E_{M_1}(e_2 S_B(b^*))
$$
which means that $\eps_B^t(b^*) = \eps_B^t(S_B(b))^*$.
Using Propositions~\ref{axiom A2}, \ref{properties of S}(iv),
and Corollary~\ref{s otimes id e} we get
$S_B(b\1) \eps_B^t(b\2) = S_B(b)H$, from where 
$HS_B(b^*) = \eps_B^t(b\2)^*S(b\1^*)$. Replacing
$S_B(b^*)$ by $b$, we get the result.
\end{proof}
\medskip

Let $\{ s_{jk}^{\alpha} \}$ be a basis consisting of matrix units of $A$
and $\{ v_{jk}^{\alpha} \}$ be a basis of comatrix units of $B$ 
dual to each other, i.e.,
$$
\la v_{jk}^{\alpha},\, s_{pq}^{\beta} \ra =
\delta_{\alpha\beta}\, \delta_{jp}\,\delta_{kq}.
$$
We have $\Delta_B(v_{jk}^{\alpha}) = \sum_l\, v_{jl}^{\alpha} \otimes 
v_{lk}^{\alpha}$ and  $\eps_B(v_{jk}^{\alpha}) =\delta_{jk}.$

\begin{lemma}
\label{v-s relations}
Let $/\alpha/ = \tau(s_{kk}^{\alpha})$.
The following identities hold true :
\begin{enumerate}
\item[(i)] 
$ E_{M_1}(e_2 e_1 v_{jk}^{\alpha}) =
d^{-1}\lambda^2 /\alpha/^{-1} s_{kj}^{\alpha}, $
\item[(ii)]
$ E_{M_1}( v_{jk}^{\alpha} e_1 e_2) =
d^{-1}\lambda^2 /\alpha/^{-1} S_A(s_{kj}^{\alpha}), $
\item[(iii)]
$ \lambda^{-1}  E_{M^\prime}(  S_A(s_{pq}^{\beta}) v_{ij}^{\alpha} e_1)
= \delta_{\alpha\beta} \delta_{ip}  v_{qj}^{\alpha}, $
\item[(iv)]
$ S_B(v_{jk}^{\alpha}) = (v_{kj}^{\alpha})^*$. 
\end{enumerate}
\end{lemma}
\begin{proof}
(i) We can directly compute :
\begin{eqnarray*}
d\lambda^{-2} /\alpha/ \tau(s_{pq}^{\beta} E_{M_1}( e_2 e_1 v_{jk}^{\alpha}) )
&=& /\alpha/ \la s_{pq}^{\beta},\,v_{jk}^{\alpha}\ra \\
&=& /\alpha/ \delta_{\alpha\beta}\, \delta_{jp}\,\delta_{kq} \\
&=& \tau( s_{kj}^{\alpha}  s_{pq}^{\beta}),
\end{eqnarray*}
therefore, we have $ E_{M_1}(e_2 e_1 v_{jk}^{\alpha}) =
d^{-1}\lambda^2 /\alpha/^{-1} s_{kj}^{\alpha}$ by the faithfulness of $\tau$.
\newline 
(ii) Similarly to (i), we compute
\begin{eqnarray*}
d\lambda^{-2} /\alpha/ \tau(E_{M_1}(v_{jk}^{\alpha} e_1 e_2)
S_A(s_{pq}^{\beta})  )
&=& /\alpha/ \la s_{pq}^{\beta},\,v_{jk}^{\alpha} \ra \\
&=& /\alpha/ \delta_{\alpha\beta}\, \delta_{jp}\,\delta_{kq} \\
&=& \tau( S_A(s_{kj}^{\alpha})  S_A(s_{pq}^{\beta}) ),
\end{eqnarray*}
and since $\tau\circ S_A = \tau$, the result follows.
\newline
(iii) Using Remark~\ref{even more}, we have
\begin{eqnarray*}
\la s_{rt}^{\gamma},\, \lambda^{-1} 
E_{M^\prime} ( S_A(s_{pq}^{\beta}) v_{ij}^{\alpha} e_1) \ra
&=&  \la s_{rt}^{\gamma},\, \lambda^{-1} 
E_{M^\prime} ( e_1 v_{ij}^{\alpha} s_{pq}^{\beta}) \ra \\
&=& \la s_{pq}^{\beta} s_{rt}^{\gamma},\, v_{ij}^{\alpha} \ra \\
&=& \delta_{\alpha\gamma}\,\delta_{qr}\,\delta_{ip}\,
\delta_{\alpha\beta}\,\delta_{tj} \\
&=&  \delta_{\alpha\beta} \delta_{ip}\, 
     \la s_{rt}^{\gamma},\, v_{qj}^{\alpha} \ra.
\end{eqnarray*}
\newline
(iv) Using part (i),  we have
\begin{eqnarray*}
 E_{M_1}( (v_{kj}^{\alpha})^* e_1 e_2) 
&=& E_{M_1}( e_2 e_1 v_{kj}^{\alpha})^*  
 = d^{-1}\lambda^2 /\alpha/^{-1} s_{kj}^{\alpha} \\
&=&  E_{M_1}( e_2 e_1 v_{jk}^{\alpha}) 
 =   E_{M_1}(S_B(v_{jk}^{\alpha})e_1 e_2  ),
\end{eqnarray*}
and the result follows from the injectivity of the map
$b\mapsto  E_{M_1}(b e_1 e_2)$.
\end{proof}

\begin{corollary}
\label{*-preserving}
$\Delta_B(b^*) = \Delta_B(b)^{*\otimes *}$, i.e., $\Delta_B$ is a
$*$-preserving map.
\end{corollary}
\begin{proof}
Using Lemmas~\ref{v-s relations}(iv) and Lemmas~\ref{properties of S}(iv),
we have
\begin{eqnarray*}
\Delta_B((v_{jk}^{\alpha})^*) 
&=& \Delta_B( S_B(v_{kj}^{\alpha}) ) 
 = \Sigma_i\, S_B(v_{ij}^{\alpha}) \otimes S_B(v_{ki}^{\alpha})  \\
&=& \Sigma_i\, (v_{ji}^{\alpha})^* \otimes (v_{ik}^{\alpha})^* 
 =  \Delta_B(v_{jk}^{\alpha})^{*\otimes *}.
\end{eqnarray*}
\end{proof} 

\begin{proposition}
\label{relations for vjk}
$ v_{ij}^{\alpha} e_1 = \lambda^{-1} \sum_k\, 
E_{M_1}(v_{ik}^{\alpha} e_1 e_2) H^{-1} v_{kj}^{\alpha}$.
\end{proposition}
\begin{proof}
By Lemma~\ref{v-s relations}(ii), all we need to show is
$$
v_{ij}^{\alpha} e_1 = d^{-1}\lambda /\alpha/^{-1}
\sum_k\,S_A(s_{ki}^{\alpha}) H^{-1} v_{kj}^{\alpha}.
$$
Since $N^\prime\cap M_2$ is spanned by the elements of the form
$ v_{rt}^{\gamma}  S_A(s_{pq}^{\beta})$, it suffices to verify that 
$$
\tau( v_{rt}^{\gamma}  S_A(s_{pq}^{\beta})  v_{ij}^{\alpha}e_1 )
=  d^{-1}\lambda /\alpha/^{-1} \sum_k\,
\tau( v_{rt}^{\gamma}  S_A(s_{pq}^{\beta}) S_A(s_{ki}^{\alpha}) 
   H^{-1} v_{kj}^{\alpha} ),
$$
or, equivalently,
$$
 E_{M^\prime}(  S_A(s_{pq}^{\beta}) v_{ij}^{\alpha}  e_1 )
= \delta_{\alpha\beta} \delta_{ip}  \lambda d^{-1} /\alpha/^{-1}
 \sum_k\, E_{M^\prime}(  S_A(s_{kq}^{\beta}) H^{-1} v_{kj}^{\alpha}  ). 
$$
Using Lemma~\ref{v-s relations}(iii), we can reduce the proof
to the verification of the relation
$$
 v_{qj}^{\alpha} =  d^{-1}  /\alpha/^{-1}
 \sum_k\,  E_{M^\prime}(  S_A(s_{kq}^{\alpha}) )  H^{-1} v_{kj}^{\alpha}.
$$
By Lemma~\ref{v-s relations}(ii),
$$ 
E_{M^\prime}(  S_A(s_{kq}^{\alpha}) ) =
d \lambda^{-2} /\alpha / E_{M^\prime}\circ E_{M_1} ( v_{qk}^{\alpha} e_1 e_2)
= d \lambda^{-1} /\alpha / E_{M_1}( v_{qk}^{\alpha} e_2), 
$$
therefore the previous relation is equivalent to
$$
v_{qj}^{\alpha} =  \lambda^{-1}             
\sum_k\,  E_{M_1}( v_{qk}^{\alpha} e_2)   H^{-1} v_{kj}^{\alpha}.
$$
Since $H\in Z(B_t)$, this is precisely Proposition~\ref{my pain}
with $b=v_{qj}^{\alpha}$, so the proof is complete.
\end{proof}

\begin{corollary}
\label{xa}
$ bx = \lambda^{-1} E_{M_1}(b\1 x e_2) H^{-1} b\2$ for all
$b\in B$ and $x\in M_1$.
\end{corollary}
\begin{proof}
Proposition~\ref{relations for vjk} implies that
$ b e_1  = \lambda^{-1} E_{M_1}(b\1 e_1 e_2) H^{-1} b\2$
for all $b\in B$. As in Remark~\ref{even more},
any $x\in M_1$ can be written as a finite sum
$x = \sum x_i e_2 y_i$ with $x_i,y_i \in M  \subset B^\prime$,
therefore, we have
\begin{eqnarray*}
bx 
&=&  \sum x_i b e_1 y_i 
 =  \sum x_i \lambda^{-1} E_{M_1}(b\1 e_1 e_2) H^{-1} b\2 y_i \\
&=& \lambda^{-1}  E_{M_1}(b\1 \sum x_i e_1  y_i e_2) H^{-1} b\2 \\
&=& \lambda^{-1}  E_{M_1}(b\1 x e_2) H^{-1} b\2.
\end{eqnarray*}
\end{proof}

\begin{proposition}
\label{the heart}
For all $x, y  \in M_1$, 
$$
E_{M_1}(bxy e_2) = \lambda^{-1} E_{M_1}(b\1 x e_2) H^{-1} E_{M_1}(b\2 y e_2).
$$
\end{proposition}
\begin{proof}
Multiplying the formula from Corollary~\ref{xa} on the right
by $ye_2t$ with $y,t\in M_1$ and taking $\tau$ from both sides we get
$$
\tau ( bxy e_2 t) = 
\lambda^{-1}\tau (E_{M_1}(b\1x e_2) H^{-1} b\2 y e_2 t),
$$
for all $t\in M_1$, from where the result follows.
\end{proof}

\begin{proposition}
\label{multiplicativity}
$\Delta_B(bc) =  \Delta_B(b) (1\otimes H^{-1}) \Delta_B(c),$
for all $b,\,c\in B.$ 
\end{proposition}
\begin{proof}
By Lemma~\ref{trick} and Proposition~\ref{the heart} 
we have for all $a_1,a_2 \in A$ : 
\begin{eqnarray*}
\lefteqn { \la a_1a_2,\, bc \ra = } \\ 
&=& \la \lambda^{-1}E_{M_1}(ca_1a_2e_2),\, b \ra \\
&=& \la \lambda^{-2}E_{M_1}(c\1 a_1e_2)H^{-1}E_{M_1}(c\2 a_2e_2),\, b \ra \\  
&=& \la \lambda^{-1}E_{M_1}(c\1 a_1e_2), b\1 \ra
    \la \lambda^{-1}E_{M_1}(H^{-1}c\2 a_2 e_2),\, b\2 \ra    \\
&=& \la a_1,\, b\1c\1 \ra \la a_2,\, b\2 H^{-1}c\2 \ra, 
\end{eqnarray*}
from where $\Delta_B(bc) =  b\1c\1 \otimes b\2 H^{-1}c\2$ which
is the result.
\end{proof}

\begin{proposition}
\label{axiom A4}
$ b\1 S_B( b\2 H^{-1}) = \eps_B^t(b)$.
\end{proposition}
\begin{proof}
Using Corollary~\ref{*-preserving}, Proposition~\ref{the heart}
and Proposition~\ref{counital maps} we have
\begin{eqnarray*}
\la a,\, b\1 S_B( b\2 H^{-1}) \ra 
&=&  d\lambda^{-3} \tau (E_{M_1}(S_B( b\2 H^{-1}) a e_2) e_2 e_1 b\1) \\
&=&  d\lambda^{-3} \tau (E_{M_1}(e_2 a b\2 H^{-1}) e_2 e_1 b\1) \\
&=& d\lambda^{-3} \tau (E_{M_1}(e_2 a b\2 H^{-1})E_{M_1}(e_2 e_1 b\1) ) \\
&=& d\lambda^{-2} \tau (E_{M_1}(e_2 a e_1 b)) \\
&=& \la a,\, \eps_B^t(b) \ra. \\ 
\end{eqnarray*}
\end{proof}

The next Corollary summarizes the properties of $\Delta_B$, $\eps_B$ and $S_B$.

\begin{corollary}
\label{bad properties}
$(\Delta_B,\,\eps_B)$ defines a coalgebra structure on $B$ such that
$$
\Delta_B(bc) =  \Delta_B(b)(1\otimes H^{-1})\Delta_B(c)
\qquad
\Delta_B(b^*) = \Delta_B(b)^{*\otimes *},
$$
the map $\eps_B^t$, defined by $\eps_B^t(b) =\eps_B(1\1b)1\2$,
satisfies the relations
$$
b\1 \otimes \eps_B^t(b\2) = 1\1b \otimes 1\2, \qquad
b\eps_B^t(c) = \eps_B(b\1 c)b\2,
$$
and there is a $*$-preserving anti-algebra and 
anti-coalgebra involution $S_B$ such that
$$
b\1 S_B( b\2H^{-1}) = \eps_B^t(b),
$$
for all $b,c\in B$.
\end{corollary}

\begin{theorem}
\label{integer case}
The following conditions are equivalent :
\begin{enumerate}
\item[(i)] $(B,\, \Delta_B,\,\eps_B,\, S_B)$ is a weak Kac algebra
with the Haar projection $e_2$ and the normalized Haar trace 
$\phi(b)= d\tau(b),\, b\in B$,
\item[(ii)] $H =1$.
%
\end{enumerate}
Moreover, if these conditions are satisfied, then $\lambda^{-1}$ is an integer.
\end{theorem}
\begin{proof}
(i)$\Rightarrow$(ii).
If $\Delta_B$ is an algebra homomorphism, then we must have 
$\Delta_B(1) = \Delta_B(1)(1\otimes H^{-1})$, and applying $(\eps_B\otimes\id)$
we get $H^{-1} =1$. 
\medskip

(ii)$\Rightarrow$(i). Clearly, if $H=1$, 
then $(B,\, \Delta_B,\,\eps_B,\, S_B)$ 
is a weak Kac algebra. For all $b\in B$ we have, by 
Proposition~\ref{counital maps} :
$$
b e_2 =\lambda^{-1} E_{M_1}(b e_2)e_2  = \eps_B^t(b)e_2,
$$
and we easily get  $S_B(e_2)=e_2$ and $\eps_B^t(e_2) =1$,
so $e_2$ is the Haar projection in $B$.

Next, since $\tau(b) =d^{-1}\la e_1,\, b\ra$, we have 
by Proposition~\ref{counital maps} :
\begin{eqnarray*}
\la a,\, \eps_B^t(b\1) \tau(b\2) \ra
&=& d^{-1} \la a ,\, \eps_B^t(b\1)\ra  \la e_1,\, b\2 \ra \\
&=& d^{-1} \la \lambda^{-1}E_{M}(ae_1) ,\, b\1\ra \la e_1,\, b\2 \ra \\
&=& d^{-1} \la \lambda^{-1}E_{M}(ae_1)e_1 ,\, b \ra \\
&=& d^{-1} \la  ae_1,\, b\ra = \la a,\, b\1 \tau(b\2) \ra,
\end{eqnarray*}
from where we get $\eps_B^t(b\1)\phi(b\2) = b\1 \phi(b\2)$. Also,
$\tau(S_B(b)) = \tau(b)$ and $\tau\circ\eps_B^t(b)
= \lambda^{-1} \tau(E_{M_1}(b e_2)) = d^{-1}\eps_B(b)$,
therefore $\phi\circ S_B = \phi$ and $\phi\circ\eps_B^t =\eps_B$.
Thus, $\phi$ is the normalized Haar trace.
\medskip

If $H=1$, then the  `trace vector' of the restriction of $\tau$ on 
$B_t$ is given by $\vec{\tau} =\frac{1}{d}(m_1,\,m_2,\dots)$,
so the componetnts of $\vec{\tau}$ are rational numbers.
Let  $\Lambda$ be the inclusion matrix of $B_t \subset B$, then
$$
\Lambda\Lambda^t \vec{\tau} = \lambda^{-1} \vec{\tau}.
$$
Since all entries of $\Lambda\Lambda^t$ and $\vec{\tau}$ are rational,
$\lambda^{-1}$ must be rational. On the other hand, $\lambda^{-1}$ is
an algebraic integer as an eigenvalue of the integer matrix
$\Lambda\Lambda^t$. Therefore, $\lambda^{-1}$ is integer.
\end{proof}

\begin{proposition}
\label{square free}
If $N\subset M$ is a depth 2 inclusion of II${}_1~$factors such that
$[M:N]$ is a square free integer (i.e., $[M:N]$ is an integer which
has no divisors of the form $n^2,\, n > 1$), then 
$N^\prime\cap M=\Bbb{C}$, and there is a (canonical) minimal action
of a Kac algebra $B$ on $M_1$ such that $M_2 \cong  M_1\rtimes B$ 
and $M=M_1^B$.
\end{proposition}
\begin{proof}  
It suffices to show that  $N\subset M$ is irreducible, since the rest 
follows from \cite{S}. Let $q$ be a
minimal projection in $M^\prime\cap M_1$, then the 
reduced inclusion $qM \subset qM_1q$ is of finite depth \cite{Bisch}.
Since any finite depth inclusion is extremal (see, e.g., \cite{Popa1},
1.3.6) we have 
$$
[qM_1q:qM] = \tau(q)^2 [M_1:M] =  \tau(q)^2 [M:N],
$$
by (\cite{PP}, Corollary 4.5).

We claim that $\tau(q)$ is a rational number. Indeed, it is well-known
that the Perron-Frobenius eigenspace of the non-negative matrix
$\Lambda\Lambda^t$ is 1-dimensional \cite{Gant}. Letting one of the
components of a corresponding eigenvector $\tau$ to be equal to
$1$, one can recover the rest of components from the system of linear 
equations with integer coefficients. Thus, we have that all components of 
$\vec{\tau}$ are rational; clearly, the normalization condition
$\tau(1)=1$ does not change this property.

Therefore, the index $[qM_1q:qM]$ is a rational number. On the other hand,
it must be an algebraic integer, since the depth is finite.
Therefore, $[qM_1q:qM]$ is an integer. Since $[M:N]$ is square free,
we must have $\tau(q)=1$, which means that 
$M^\prime \cap M_1$ and $N^\prime\cap M$ are 1-dimensional.
\end{proof} 

\begin{corollary}
If $N\subset M$ is a depth 2 inclusion of II${}_1~$factors such that
$[M:N]=p$ is prime, then $N^\prime\cap M=\Bbb{C}$, and there is an
outer action of the cyclic group $G=Z/pZ$ on $M_1$ such that 
$M_2 \cong  M_1\rtimes G$ and $M=M_1^G$.
\end{corollary}
\begin{proof}
By Proposition~\ref{square free}, $B$ must be a Kac algebra of prime
dimension $p$. But it is known that any such an algebra is a group algebra
of the cyclic group $G=Z/pZ$ \cite{Kac}.
\end{proof}

\end{section}


\begin{section}
{Weak $C^*$-Hopf algebra structure on $M'\cap M_2$ 
(the general  case) }

When $H\neq 1$,  $(B,\, \Delta_B,\,\eps_B,\, S_B)$ is no longer 
a weak Kac algebra (for instance, $\Delta_B$ is not a homomorphism).
However, it is possible to deform the structure maps in such a way that
$A$ becomes a weak $C^*$-Hopf algebra.

\begin{definition}
\label{correct maps}
Let us define the following operations on $B$ :
$$
\begin{array}{rrrl}
involution      & \star : B\to B &:& b^\star = S_B(H)^{-1}b^*S_B(H), \\
comultiplication& \tDelta : B\to B\otimes B &:&
\tDelta(b) = (1\otimes H^{-1})\Delta_B(b)~i.e.,  \\
 &  &  & b\tI \otimes b\tII = b\1 \otimes H^{-1} b\2 \\
counit & \teps : B\to \Bbb{C} &:& \teps(b) =\eps_B(Hb), \\
antipode & \tS : B\to B &:&  \tS(b) = S_B(HbH^{-1}).
\end{array}
$$
\end{definition}
\medskip

Clearly, $\star$ defines a $C^*$-algebra structure on $B$ (we will
still denote this new $C^*$-algebra by $B$). Our goal is to show that
$(B,\, \tDelta,\,\teps,\, \tS)$ is a weak $C^*$-Hopf algebra. The proof
of this fact consists of a verification of all the axioms from Section 2. 
We will need the following technical lemma.

\begin{lemma}
\label{a and z}
For all $b\in B$ and $z\in B_t$ we have
\begin{enumerate}
\item[(i)] $\eps_B^t(zb) = z \eps_B^t(b)$,
\item[(ii)] $b\1 z\otimes b\2 =(bz)\1 \otimes (bz)\2$,
\item[(iii)] $ b\1S_B(z) \otimes b\2 = b\1 \otimes b\2 z$,
\end{enumerate}
\end{lemma}
\begin{proof}
Part (i) is clear from Proposition~\ref{counital maps}.
Next, recall that $B_t =A \cap B$, and compute
$$
\la a_1,\, b\1 z\ra \la a_2,\, b\2\ra = \la za_1a_2,\, b \ra 
=\la a_1a_2,\, bz \ra, \quad a_1,a_2\in A,
$$
which gives (ii). Finally, using the properties of $S_B$ we have
\begin{eqnarray*}
\la a_1,\, b\1S_B(z)\ra \la a_2,\, b\2 \ra
&=& \la a_1,\, S_B(z S_B(b\1))\ra \la a_2,\, b\2 \ra \\ 
&=& \overline{ \la a_1^*,\, S_B(b\1^*) z^* \ra} \la a_2,\, b\2 \ra \\
&=& \overline{ \la (a_1z)^*,\, S_B(b\1^*) \ra} \la a_2,\, b\2 \ra \\
&=& \la a_1z,\, b\1 \ra \la a_2,\, b\2 \ra \\ 
&=& \la a_1za_2,\, b \ra \\
&=& \la a_1,\, b\1 \ra \la a_2,\, b\2z \ra,
\end{eqnarray*}
from where (iii) follows. 
\end{proof}

\begin{proposition}
\label{coalgebra axiom}
$(B,\, \tDelta,\,\teps)$ is a coalgebra.
\end{proposition}
\begin{proof}
Let us check the coassociativity of $\tDelta$. Using Lemma~\ref{a and z}
and the fact that $H\in B_t$ we compute for all $b\in B$ :
\begin{eqnarray*}
(\tDelta\otimes \id)\tDelta(b)
&=& \tDelta(b\1) \otimes H^{-1}b\2 
     =  b\1  \otimes H^{-1}b\2 \otimes H^{-1} b\3 \\
&=& b\1 \otimes (H^{-1}b\2)\1 \otimes H^{-1}(H^{-1}b\2)\2 \\
&=& b\1 \otimes \tDelta(H^{-1}b\2) = (\id \otimes \tDelta)\tDelta(b).  
\end{eqnarray*}
Next, we check the counit axioms :
\begin{eqnarray*}
(\teps\otimes\id)\tDelta(b) &=& \eps(H b\1)H^{-1}b\2 
= \eps((Hb)\1)H^{-1}(Hb)\2 = b, \\
(\id\otimes\teps)\tDelta(b) &=& b\1 \eps(HH^{-1}b\2) =b. 
\end{eqnarray*}
\end{proof}

\begin{proposition}
\label{comultiplication axioms}
$\tDelta$ is a $\star$-homomorphism.
\end{proposition}
\begin{proof}
Using the properties of $\Delta_B$ from Corollary~\ref{bad properties}
and Lemma~\ref{a and z} we have:
\begin{eqnarray*}
\tDelta(bc)
&=& (1 \otimes H^{-1}) \Delta_B(bc) \\
&=& (1 \otimes H^{-1}) \Delta_B(b) (1 \otimes H^{-1}) \Delta_B(c)
= \tDelta(b)\tDelta(c), \\
\tDelta(b^\star)
&=& \tDelta(S_B(H)^{-1} b^*S_B(H)) \\
&=& (S_B(H)^{-1}b^*S_B(H))\1 \otimes H^{-1} (S_B(H)^{-1}b^*S_B(H))\2 \\
&=& S_B(H)^{-1}b\1^* \otimes S_B(H)^{-1} b\2^* S_B(H) 
 = (S_B(H)^{-1} b\1)^\star \otimes b\2^\star \\
&=& b\1^\star \otimes (H^{-1}b\2)^\star = \tDelta(b)^{\star\otimes\star}.
\end{eqnarray*}
\end{proof}

\begin{proposition}
\label{counital axioms}
Let $\teps^t(b)= \teps(1\tI b)1\tII$. Then, for all $b,\,c\in B$ :
$$
b\teps^t(c) = \teps(b\tI c)b\tII, \qquad
b\tI \otimes \teps^t(b\tII) =  1\tI b\otimes 1\tII,
$$
\end{proposition}
\begin{proof}
First, we compute, using Lemma~\ref{a and z} and Proposition~\ref{axiom A2} :
$$
\teps^t(b) = \eps(H1\1b)H^{-1}1\2 = \eps(H\1b)H\2H^{-1} =
H\eps_B^t(b)H^{-1} = \eps_B^t(b).
$$ 
Using this relation, Lemma~\ref{a and z}, and properties of 
$\eps_B^t$ from Corollary~\ref{bad properties} we have 
\begin{eqnarray*}
b\tI \otimes \teps^t(b\tII) 
&=& b\1 \otimes \eps_B^t(H^{-1}b\2) = b\1 \otimes H^{-1} \eps_B^t(b\2) \\
&=& 1\1 b \otimes H^{-1} 1\2 = 1\tI b \otimes 1\tII, \\
b \teps^t(c) 
&=& b \eps_B^t(c) = H^{-1}(Hb)\eps_B^t(c) \\
&=& H^{-1} \eps_B((Hb)\1c)(Hb)\2 = \eps_B(Hb\1c) H^{-1} b\2 \\
&=& \teps(b\tI c) b\tII.
\end{eqnarray*}
\end{proof}
   
\begin{proposition}
\label{antipode axioms}
$\tS$ is a linear anti-multiplicative and 
anti-comultiplicative map such that 
$$
b\tI \tS(b\tII)  = \teps^t(b).
$$
Moreover, $(\tS\circ\star)^2=\id$ and 
$\tS^2(b) = GbG^{-1}$, where $G=\tS(H)^{-1}H$.  
\end{proposition}  
\begin{proof}
Using Corollary~\ref{bad properties}, Lemma~\ref{a and z} and
 definitions of $\tS$ and $\star$, we have :
\begin{eqnarray*}
\tS(bc) 
&=&  S_B(H bc H^{-1}) = S_B(HcH^{-1})S_B(HbH^{-1})= \tS(c)\tS(b), \\
\tS(b)\tII \otimes \tS(b)\tI
&=& H^{-1} S_B(HbH^{-1})\2 \otimes S_B(HbH^{-1})\1 \\
&=& S_B(HbH^{-1})\2 \otimes S_B(H^{-1})S_B(HbH^{-1})\1 \\
&=& S_B((HbH^{-1})\1) \otimes S_B(H^{-1}) S_B((HbH^{-1})\2) \\
&=& S_B(H b\1H^{-1}) \otimes S_B(b\2H^{-1}) \\       
&=& \tS(b\1)  \otimes \tS(H^{-1}b\2) = \tS(b\tI) \otimes \tS(b\tII), \\
b\tI \tS(b\tII) 
&=& b\1 S_B(b\2 H^{-1}) =\eps_B^t(b) = \teps^t(b), 
\end{eqnarray*}
from where the first part of Proposition follows. 
Next, we can compute
\begin{eqnarray*}
\tS(b^\star) 
&=& S_B(H b^\star H^{-1}) = S_B(H S_B(H)^{-1} b^* S_B(H)H^{-1}) \\
&=& S_B(H)^{-1} H S_B(b^*) H^{-1} S_B(H), \\
\tS(\tS(b^\star)^\star)
&=& \tS((S_B(H)^{-1} H S_B(b^*) H^{-1} S_B(H))^\star) = \tS(H^{-1} S_B(b) H) \\
&=& S_B(S_B(b)) = b,
\end{eqnarray*}
therefore $(\tS\circ\star)^2=\id$. Finally, since $S_B(H)=\tS(H)$, we get
\begin{eqnarray*}
\tS^2(b) 
&=&\tS(S_B(H bH^{-1})) = S_B(H S_B(H bH^{-1}) H^{-1}) \\
&=& S_B(H)^{-1} H b  H^{-1}S_B(H) = GbG^{-1}.
\end{eqnarray*}
\end{proof}
 
Thus, we can state the main result of this section.

\begin{theorem}
\label{WHA}
$(B,\, \tDelta,\,\teps,\, \tS)$ is a weak $C^*$-Hopf algebra
with the Haar projection $e_2H$ and  normalized Haar functional
$\tilde\phi(b) = \phi (H\tS(H)b) = d\tau(\tS(H)Hb)$ 
(cf.\ Theorem~\ref{integer case} ).
\end{theorem}
\begin{proof}
It follows from Propositions~\ref{coalgebra axiom}
-- \ref{antipode axioms} that $(B,\, \tDelta,\,\teps,\, \tS)$ 
is a weak $C^*$-Hopf algebra. 

The properties of $e_2$ established in  Theorem~\ref{integer case} 
and Proposition~\ref{counital maps} give
\begin{eqnarray*} 
be_2H 
&=& \eps_B^t(b) e_2 H = \teps^t(b)e_2H, \\
\teps^t(e_2H) 
&=& \eps_B^t(e_2H) =\lambda^{-1}E_{M_1}(e_2He_2) \\
&=& \lambda^{-1}E_{M_1}(E_M(H)e_2) = E_M(H) =1,
\end{eqnarray*}
since $\teps^t =\eps_B^t$ by Proposition~\ref{counital axioms},
and $E_M(H) =\tau(H)1 =1$ by Corollary~\ref{s otimes id e}.

Using Lemma~\ref{a and z}(ii), and taking into account that 
$\tS^2\vert_{B_t} = \id_{B_t}$ (Proposition~\ref{antipode axioms}), 
we compute for all  $b\in B,\, z\in B_t$ :
$$
\teps^t(\tS(z)) = \tS(z)\tI \tS(\tS(z)\tII) =1\tI \tS(\tS(z) 1\tII) 
 = \tS^2(z) =z, 
$$
therefore $e_2 \tS(H) = e_2 \teps^t(\tS(H))^* = e_2 H$. 
Since $S_B(e_2) =e_2$ and $S_B(H) =\tS(H)$, using the above relation,
we get
$$
\tS(e_2H) = \tS(H)^{-1} S_B(e_2H) \tS(H) = e_2 \tS(H) = e_2 H.
$$
Thus $\tS(e_2H) =e_2 H$. Also we have :
\begin{eqnarray*}
(e_2H)^2 &=& E_M(H)e_2 H = e_2H, \\
(e_2H)^\star &=& \tS(H)^{-1} H e_2 \tS(H) = e_2 \tS(H) = e_2 H.
\end{eqnarray*}
Therefore, $e_2H$ is an $\tS$-invariant projection. 
This proves that $e_2$ is the Haar projection of $B$.

Next, using Lemma~\ref{a and z} and the properties of the trace 
$\phi$ from the proof of Theorem~\ref{integer case} we have
\begin{eqnarray*}
\teps^t(b\tI)\tilde\phi(b\tII)
&=& \eps_B^t(b\1)\tilde\phi(H^{-1}b\2) 
 = \eps_B^t(b\1)\phi(\tS(H)b\2) \\
&=& \eps_B^t( (b\tS(H))\1) \phi((b\tS(H))\2)  
 = (b\tS(H)))\1  \phi((b\tS(H))\2)  \\
&=& b\1 \phi(\tS(H)b\2) 
 = b\1 \phi(H\tS(H)H^{-1} b\2) \\
&=& b\1 \tilde\phi(H^{-1} b\2) 
 = b\tI \tilde\phi(b\tII), \\ 
\tilde\phi(\tS(b)) 
&=& \phi (H\tS(H)\tS(H)^{-1} S_B(b) \tS(H)) = \phi(H\tS(H) S_B(b)) \\
&=& \phi(b \tS(H) H) = \tilde\phi(b), \\
\tilde\phi(\teps^t(b)) 
&=& \phi(\tS(H) H\eps_B^t(b)) = \tau(\tS(H)) \phi(\eps_B^t(Hb)) 
=\eps_B(Hb) =\teps(b),
\end{eqnarray*}
therefore, $\tilde\phi$ is the normalized Haar functional on $B$.
\end{proof}

\begin{remark}
\begin{enumerate}
\item[(i)]
The non-degenerate duality $\la,\,\ra$ induces on $A=N^\prime \cap M_1$ 
the  structure of the weak $C^*$-Hopf algebra dual to $B$.
\item[(ii)]
The weak $C^*$-Hopf algebra $B$ is biconnected, since the inclusion
$B_t=M^\prime \cap M_1 \subset B =M^\prime \cap M_2$  is connected 
(\cite{GHJ}, 4.6.3) and $B_t \cap B_s =(M^\prime \cap M_1)\cap 
(M_1^\prime \cap M_2)= \Bbb{C}$. Thus, only biconnected weak Hopf $C^*$-algebras 
arise as symmetries of finite index depth 2 inclusions of II${}_1~$factors.
\item[(iii)] 
If $\lambda^{-1}$ is not integer, then $\tS$ has infinite order.
Indeed, the canonical element $G$ implementing the square of the
antipode in Proposition~\ref{antipode axioms} is positive, so if
$\tS^{2n} =\id$ for some $n$, then $G^n$ belongs to $Z(B)$, the center of $B$.
Taking the $n$-th root, we get that $G\in Z(B)$, which means that $S^2=\id$,
and $B$ is a weak Kac algebra, which is in contradiction with 
Theorem~\ref{integer case}.
\end{enumerate}
\end{remark}

\end{section}


\begin{section}
{Action of $B$ on $M_1$.}

Note that in terms of $\tDelta$, Proposition~\ref{the heart} means
that
$$
E_{M_1}( b xy e_2) =\lambda^{-1} E_{M_1}( b\tI x e_2) E_{M_1}( b\tII y e_2),
$$
for all $b\in B$ and $x,y\in M_1$. This suggests the following definition 
of the action of $B$ on $M_1$.

\begin{proposition}
\label{the action}
The map $\lact : B \otimes M_1 \to M_1$ :
$$
b \lact x = \lambda^{-1} E_{M_1}(b x e_2)
$$
defines a left action of $B$ on $M_1$ (cf.\ \cite{S}, Proposition 17).
\end{proposition}
\begin{proof}
Clearly, the above map defines a left $B$-module structure
on $M_1$, since $1\lact x =x$ and
$$
b \lact (c \lact x) 
= \lambda^{-2}  E_{M_1}(b  E_{M_1}(c x e_2)  e_2) 
= \lambda^{-1}  E_{M_1}( bc x e_2 ) = (bc) \lact x.
$$
Next, using Proposition~\ref{the heart} we get
\begin{eqnarray*}
b \lact xy 
&=& \lambda^{-1} E_{M_1}( b xy e_2)  
= \lambda^{-2} E_{M_1}( b\tI x e_2) E_{M_1}( b\tII y e_2) \\
&=& ( b\tI \lact x)( b\tII \lact y).
\end{eqnarray*}
By Remark~\ref{even more} and properties of $S_B$ we also get
\begin{eqnarray*}
\tS(b)^\star \lact x^*
&=& \lambda^{-1} E_{M_1}( \tS(b)^\star x^* e_2 ) \\
&=& \lambda^{-1} E_{M_1}( S_B(H)^{-1} S_B(HbH^{-1})^* S_B(H) x^* e_2 ) \\
&=& \lambda^{-1} E_{M_1}( S_B(b^*)  x^* e_2 ) \\
&=&  \lambda^{-1} E_{M_1}( e_2  x^* b^* ) 
 =  \lambda^{-1} E_{M_1}( bx e_2)^* = (b \lact x)^*.
\end{eqnarray*}
Finally,
$$
b \lact 1 = \lambda^{-1} E_{M_1}( b e_2) = 
\lambda^{-1} E_{M_1}( \lambda^{-1} E_{M_1}(b e_2)e_2)
= \teps^t(b) \lact 1,
$$ 
and $b\lact 1 =0$  iff $\teps^t(b) = \lambda^{-1} E_{M_1}(b e_2)= 0$.
\end{proof}

\begin{proposition}
\label{fixed point}
$M_1^B = M$, i.e. $M$ is the fixed point subalgebra of $M_1$.
\end{proposition}
\begin{proof}
If $x\in M_1$ is such that $b\lact x = \teps^t(b)\lact x$ for all
$b \in B$, then $ E_{M_1}(bxe_2) = E_{M_1}(\eps_B^t(b)xe_2) =
E_{M_1}(be_2) x$. Taking $b =e_2$, we get $E_{M}(x) =x$ which
means that $x \in M$. Thus, $M_1^B \subset M$.

Conversely, if $x\in M$, then $x$ commutes with $e_2$ and 
$$
b\lact x = \lambda^{-1} E_{M_1}(b e_2 x) = 
\lambda^{-1} E_{M_1}( \lambda^{-1} E_{M_1}(be_2)e_2 x)
= \eps_B^t(b) \lact x,
$$
therefore $M_1^B = M$. 
\end{proof} 

\begin{proposition}
\label{crossproduct}
The map $\theta : [x\otimes b] \mapsto x \tS(H)^\half b  \tS(H)^{-\half}$
defines a von Neumann algebra isomorphism between $M_1\rtimes B$ and $M_2$.
\end{proposition}
\begin{proof}
By definition of the action $\lact$ we have :
\begin{eqnarray*}
\theta([x(z\lact 1) \otimes b]) 
&=& x \tS(H)^\half \lambda^{-1} E_{M_1}(ze_2)b \tS(H)^{-\half} \\
&=& x \tS(H)^\half zb \tS(H)^{-\half} = \theta([x\otimes zb]),
\end{eqnarray*}
for all $x \in M_1,\, b\in B,\, z\in B_t$, so
$\theta$ is a well defined linear map from
$M_1\rtimes B = M_1 \otimes_{B_t} B$  to $M_2$.
 It is surjective since an orthonormal basis of $B=M^\prime \cap M_2$
over $B_t=M^\prime \cap M_1$ is also a basis of $M_2$ over $M_1$
(\cite{Popa1}, 2.1.3).

Let us check that $\theta$ is an involution-preserving isomorphism
of algebras. Note that from Corollary~\ref{xa} we have 
$bx = (b\tI \lact x)b\tII$. This allows us to compute, 
for all $x,y\in M_1$ and $b,c\in B$:
\begin{eqnarray*}
\theta([x\otimes b][y\otimes c])
&=& \theta([x(b\tI\lact y) \otimes b\tII c]) \\
&=& x (b\tI\lact y) \tS(H)^\half b\tII c  \tS(H)^{-\half} \\
&=& x( (\tS(H)^\half b)\tI \lact y) (\tS(H)^\half b)\tII c \tS(H)^{-\half} \\
&=& x \tS(H)^\half b y c \tS(H)^{-\half} \\
&=& x \tS(H)^\half b \tS(H)^{-\half} y \tS(H)^\half c \tS(H)^{-\half} \\
&=& \theta([x\otimes b])\theta([y\otimes c]), \\
\theta([x\otimes b]^*)
&=& (b\1^\star \lact x^*) \tS(H)^\half b\2^\star  \tS(H)^{-\half} \\
&=& (\tS(H)^\half b^\star)\tI \lact x^*) 
    (\tS(H)^\half b^\star)\tII \tS(H)^{-\half} \\
&=& \tS(H)^{\half}  b^\star \tS(H)^{-\half} x^* \\
&=& \tS(H)^{-\half}  b^* \tS(H)^{\half} x^* \\
&=& ( x \tS(H)^{\half} b \tS(H)^{-\half} )^* = \theta([x\otimes b])^*.
\end{eqnarray*}

It is known that $M_1\rtimes B$ is a II${}_1$~factor iff
$M_1^B$ is \cite{NSzW}. Now the injectivity of $\theta$ follows
from the simplicity of II${}_1$~factors (see, e.g., the
appendix of \cite{JS}). Thus, $\theta$ is a von Neumann algebra
isomorphism.
\end{proof}
  
\begin{remark}
\label{minimality}
\begin{enumerate}
\item[(i)]
The action of $B$ constructed in Proposition~\ref{the action}
is minimal, since we have $M_1^\prime \cap M_1\rtimes B = M_1^\prime \cap M_2 
=B_s$ by Proposition~\ref{crossproduct}. 
\item[(ii)]
If the inclusion $N\subset M$ is irreducible, then $B$ is a usual
Kac algebra (i.e., a Hopf $C^*$-algebra) and we recover the well-known
result proved in \cite{S}, \cite{L}, and \cite{D}.
\end{enumerate}
\end{remark}

\end{section}


\end{document}